\documentclass{article}

\usepackage{delarray,verbatim,enumerate,a4wide}
\usepackage{amsmath,amsthm,amstext,amsbsy,amssymb,amsfonts,amscd}
\usepackage[all]{xy}
\usepackage[frenchb]{babel}
\usepackage[T1]{fontenc}
\usepackage[utf8]{inputenc}

\usepackage{scalerel}

\usepackage[small]{titlesec}

\usepackage{color}
\definecolor{vert}{rgb}{0.1,0.4,0.2}
\usepackage[colorlinks=true,linkcolor=blue,citecolor=vert]{hyperref}

\usepackage{calligra}
\DeclareFontShape{T1}{calligra}{m}{n}{<->s*[0.95]callig15}{}
\DeclareMathAlphabet{\mathscr}{T1}{calligra}{m}{n}

\newtheorem{Th}{Théorème}[]

\newtheorem{Prop}[Th]{Proposition}
\newtheorem{Cor}[Th]{Corollaire}

\newtheorem{Sco}[Th]{Scolie}
\newtheorem{Def} [Th]{Définition}

\newtheorem*{CoCy}{Conjecture cyclotomique}

\def\Preuve{\noindent {\it Preuve.~}}

\def\Remarque{\smallskip\noindent {\it Remarque.~}}
\def\Remarques{\smallskip\noindent {\it Remarques.~}}

\def\Nota{\smallskip\noindent {\it Nota.~}}

\def\NN{\mathbb N}	\def\ZZ{\mathbb Z}		
\def\F2{\mathbb{F}_2}	\def\Z2{\mathbb{Z}_2}		
\def\Zl{\mathbb{Z}_\ell} 			

 		\def\P{\mathcal  P}		\def\U{\mathcal  U}	\def\F{\mathcal  F}
\def\J{\mathcal  J}  	\def\C{\mathcal  C}	\def\R{\mathcal  R}	\def\X{\mathcal  X}
\def\Dl{\mathcal  D\ell} 	\def\Pl{\mathcal  P\ell}  	\def\Cl{\mathcal  C\!\ell}	
\def\E{\mathcal  E}		\def\T{\mathcal  T}		\def\Dl{\mathcal  D\ell}	\def\D{\mathcal  D}

\def\G{\mathscr G\,}

		\def\p{{\mathfrak p}}						
		\def\l{{\mathfrak l}}				\def\A{{\mathfrak A}}

	\def\deg{\operatorname{deg}}
\def\Gal{\operatorname{Gal}}

\newcommand\scale[2]{\vstretch{#1}{\hstretch{#1}{#2}}}

\newcommand\si[1]{\scale{.8}{#1}}
\newcommand\ph{{\phantom{*}}}
\newcommand\ab{{\scale{.8}{\rm ab}}}

\makeatletter
\newcommand*\wt[2][0.2ex]{%
        \begingroup
        \mathchoice{\wt@helper{#1}{#2}{\displaystyle}{\textfont}}
                   {\wt@helper{#1}{#2}{\textstyle}{\textfont}}
                   {\wt@helper{#1}{#2}{\scriptstyle}{\scriptfont}}
                   {\wt@helper{#1}{#2}{\scriptscriptstyle}{\scriptscriptfont}}%
        \endgroup
        #2%
}
\newcommand*\wt@helper[4]{%
        \def\currentfont{\the#41}%
        \def\currentskewchar{\char\the\skewchar\currentfont}%
        \setbox\tw@\hbox{\currentfont$#2$\currentskewchar}%
        \dimen@ii\wd\tw@
        \setbox\tw@\hbox{\currentfont$#2${}\currentskewchar}%
        \advance\dimen@ii-\wd\tw@
        \rlap{\raisebox{-#1}{$\m@th#3\kern\dimen@ii\widetilde{\phantom{#2}}$}}%
}
\makeatother

\def\wE{\,\wt[0.1ex]{\!\mathcal E}}		\def\wU{\wt[0.2ex]{\mathcal U}}
\def\wJ{\,\wt[0.2ex]{\!\mathcal J}}	\def\wCl{\wt[0.1ex]{\mathcal C\!\ell}} \def\wDl{\wt[0.2ex]{\mathcal D\!\ell}}
				
		\def\wP{\,\wt[0.1ex]{\!\mathcal P}}

\begin{document}

\title{\Large\bf Généralisation d'un théorème de Greenberg}

\author{ Jean-François {\sc Jaulent} }
\date{}
\maketitle
\bigskip\bigskip

{\small
\noindent{\bf Résumé.} Nous formulons une conjecture générale sur le polynôme caractéristique des modules d'Iwasawa $S$-décomposés $T$-ramifiés au-dessus de la $\Zl$-extension cyclotomique d'un corps de nombres. Nous montrons qu'elle est en fait équivalente à la conjonction des conjectures de Leopoldt et de Gross-Kuz'min. Nous étendons ainsi un théorème de semi-simplicité de Greenberg et, au passage, un isomorphisme de Kuz'min.}

\

{\small
\noindent{\bf Abstract.} We formulate a general conjecture on the characteristic polynomials of $S$-decomposed $T$-ramified Iwasawa modules over the cyclotomic $\Zl$-extension of a number field. We show that this conjecture is equivalent to the conjunctions of the classical conjectures of Leopoldt and of Gross-Kuz'min.  We so extend  a result of semi-simplicity of Greenberg and, by the way, an isomorphism of Kuz'min.}
\bigskip\bigskip\bigskip


\section*{Introduction}
\addcontentsline{toc}{section}{Introduction}
\medskip

Soit $\ell$ un nombre premier arbitraire, $K$ un corps de nombres, $K_{\si{\infty}}$ sa $\Zl$-extension cyclotomique et  $\Lambda=\Zl[[\gamma-1]]$ l'algèbre d'Iwasawa  du groupe $\Gamma=\Gal(K_{\si{\infty}}/K)=\gamma^\Zl$.\smallskip

Un théorème de Greenberg (cf. \cite{Grb1}) affirme que, si $K$ est abélien, le polynôme caractéristique du groupe de Galois $\Gal(K_{\si{\infty}}^{\si{\rm cd}}/K_{\si{\infty}})$ de la plus grande pro-$\ell$-extension abélienne non ramifiée et $\ell$-décomposée $K_{\si{\infty}}^{\si{\rm cd}}$ de  $K_{\si{\infty}}$, regardé comme $\Lambda$-module, n'est pas divisible par $(\gamma-1)$.\smallskip

L'objet de cette note est d'introduire une conjecture générale, pour tout corps de nombres $K$ et tout nombre premier $\ell$, faisant intervenir deux ensembles finis disjoints $S$ et $T$  de places finies du corps $K$, qui généralise la situation considérée par Greenberg.
Plus précisément, nous postulons que  le groupe de Galois $\,\C^T_S=\Gal(H^T_S(K_{\si{\infty}})/K_{\si{\infty}})$ attaché à la pro-$\ell$-extension abélienne maximale $H^T_S(K_{\si{\infty}})$ de $K_{\si{\infty}}$ qui est $S$-décomposée et $T$-ramifiée (i.e. complètement décomposée aux places au-dessus de $S$ et non-ramifiée en dehors ds places au-dessus de $T$) est encore un $\Lambda$-module (non nécessairement de torsion) dont le polynôme caractéristique n'est pas divisible par $(\gamma-1)$ dès lors que la réunion $S\cup T$ contient l'ensemble $Pl_\ell$ des places au-dessus de $\ell$.\smallskip

Sans surprise, cette conjecture contient de fait les conjectures de Leopoldt et de Gross-Kuz'min (cf. e.g. \cite{J55}), qui correspondent respectivement aux cas $(S,T)=(\emptyset,Pl_\ell)$ et $(S,T)=(Pl_\ell,\emptyset)$, cette dernière situation étant précisément celle considérée par Greenberg. Le résultat principal de notre étude (Th. \ref{ThP}) est que la réciproque est vraie: si $K$ est totalement réel, ou si $K$ est un corps à conjugaison complexe extension quadratique totalement imaginaire d'un sous-corps totalement réel, la conjecture cyclotomique que nous avançons est vraie dès lors que $K$ vérifie à la fois la conjecture de Leopoldt et celle de Gross-Kuz'min. En d'autres termes, la conjecture cyclotomique est  équivalente à la conjonction des conjectures de Leopoldt et de Gross-Kuz'min. Elle est de ce fait satisfaite par les corps abéliens (ce qui généralise le résultat de Greenberg) et quelques autres, notamment les corps logarithmiquement principaux, dont on sait qu'ils vérifient les conjectures de Gross-Kuz'min  et de Leopoldt en présence des racines $2\ell$-ièmes de l'unité (cf. e.g. \cite{J17}).
\medskip

\Nota Pendant la rédaction de ce travail, nous avons eu connaissance d'une pré-publication de Lee et Seo \cite{LS} qui redémontre (sans y référer explicitement) quelques-uns des résultats de \cite{J18} et propose une conjecture voisine de la nôtre (mais {\em de facto} équivalente) dans une formulation plus compliquée en l'absence de l'hypothèse $Pl_\ell\subset S \cup T$ (cf. la remarque {\em in fine} de la présente note).
\bigskip


\newpage
\section{Le $\ell$-groupe des  $S$-classes $T$-infinitésimales d'un corps de nombres}

Rappelons succinctement quelques éléments de la Théorie $\ell$-adique du corps de classes telle qu'exposée dans \cite{J31}: Le nombre premier $\ell$ étant supposé fixé, pour chaque place $\p$ du corps de nombres $K$ nous notons $\R_{K_\p}=\varprojlim \,K_\p^\times/K_\p^{\times\ell^n}$ le $\ell$-adifié du groupe multiplicatif du complété $K_\p$ de $K$ en $\p$; par $\,\U^{\phantom{*}}_{K_\p}$ son sous-groupe unité (au sens habituel); et par $\,\wU^{\phantom{*}}_{K_\p}$ le groupe des unités logarithmiques, i.e. le  sous-groupe de normes attaché à la $\Zl$-extension cyclotomique $K_\p^c$ de $K_\p$. \smallskip

Le $\ell$-adifié du groupe des idèles du corps $K$ est le produit restreint $\J(K)=\prod_\p^{\si{\rm res}}\R_{K_\p}$ formé des familles $(x_\p)_\p$ dont presque tous les éléments sont des unités. Son sous-groupe unité est le produit $\,\U(K)=\prod_\p\U_{K_\p}$; et celui des unités logarithmiques est le produit $\,\wU(K)=\prod_\p\wU_{K_\p}$. Le sous-groupe des idèles principaux est le tensorisé $\R(K)=\Zl\otimes_\ZZ K^\times$ du groupe multiplicatif de $K$; il se plonge canoniquement dans $\J(K)$ et définit un quotient compact que la  Théorie identifie au groupe de Galois $\G(K)=\Gal(K^{\ab}/K)$ de la pro-$\ell$-extension abélienne maximale de $K$.\smallskip

Cela étant, pour tout couple $(S,T)$ d'ensembles finis disjoints de places de $K$ nous notons  $\R_S(K)=\prod_{\p\in S}\R_{K_\p}$ et $\,\J^T(K)=\prod^{\si{\rm res}}_{\p\notin T}\R_{K_\p}$; puis $\,\U_S(K)=\prod_{\p\in S}\U_{K_\p}$ et $\,\U^T(K)=\prod_{\p\notin T}\U_{K_\p}$.
Et la Théorie $\ell$-adique du corps de classes identifie le groupe de Galois $\G^T_S(K)=\Gal(H^T_S(K)/K)$ de la pro-$\ell$-extension abélienne  $S$-décomposée $T$-ramifiée maximale $H^T_S(K)$ de $K$ au quotient:\smallskip

\centerline{$\J(K)/\J_S^T(K)\R(K)$, avec $\J_S^T(K)=\R_S(K)\,\U^T(K)=\prod_{\p\in S}\R_{K_\p}\prod_{\p\notin T}\U_{K_\p}$.}\smallskip

Celui-ci s'interprète alors comme groupe de classes de diviseurs de la façon suivante:

\begin{Def}
Soient $S$ et $T$ deux ensembles finis disjoints de places d'un corps de nombres $K$.
\begin{itemize}
\item Le $\ell$-groupe des $S$-diviseurs étrangers à $T$ est le quotient $\,\D_S^T(K)=\J^T(K)/\J^T_S(K)$.
\item Son sous-groupe principal $T$-infinitésimal est l'image $\P_S^T(K)= \R^T(K)\J_S^T(K)/\J^T_S(K)$ du groupe des idèles principaux $T$-infinitésimaux
$\R^T(K)=\R(K)\cap\J^T(K)$.
\item Et nous disons que le quotient $\,\Cl^T_S(K)=\D^T_S(K)/\P^T_S(K)\simeq\G_S^T(K)$ est le $\ell$-groups des $S$-classes $T$-infinitésimales du corps $K$.
\end{itemize}
\end{Def}

Bien entendu, c'est le théorème d'approximation simultanée qui permet de représenter chaque classe de $\J(K)/\R(K)$ par un élément de $\J^T(K)$.
On tombe alors sur la description du groupe $\,\Cl^T_S(K)$ donnée dans \cite{J18} (pp. 148--150), à l'inversion près de $S$ et de $T$; ce qui permet de le traiter comme un groupe de classes d'idéaux et donc de lui appliquer en particulier les calculs de classes invariantes à la Chevalley qui conduisent à une généralisation naturelle de la suite exacte des classes ambiges dans ce nouveau contexte (cf. \cite{J18}, Th. II.2.33). En particulier, il vient:

\begin{Prop}\label{Prop}
Dans une $\ell$-extension cyclique $N/K$ de corps de nombres, le quotient du $\ell$-groupe des $S$-classes $T$-infinitésimales ambiges par le sous-groupe des classes des $S$-diviseurs ambiges étrangers à $T$ est donné par l'isomorphisme:\smallskip

\centerline{$\Cl^T_S(N)^\Gamma/cl^T_S\big(\D^T_S(N)^\Gamma\big) \simeq \big(\E^T_S(K)\cap N_{N/K}(\R(N))\big)/N_{N/K}(\E_S^T(N))$.}\smallskip

\noindent Ici $\,\E^T_S(K)\cap N_{N/K}(\R(N))$ est le pro-$\ell$-groupe des $S$-unités $T$-infinitésimales du corps $K$ qui sont normes dans $N/K$ et $N_{L/K}(\E_S^T(N))$ est le sous-groupe des normes de $S$-unités $T$-infinitésimales.
\end{Prop}

\Preuve Pour la commodité du lecteur, plutôt que de renvoyer aux arguments cohomologiques développés dans \cite{J18}, expliquons simplement comment ce résultat s'explicite dans le cas cyclique: tout comme dans la preuve historique de Chevalley pour les classes d'idéaux prenons un générateur arbitraire $\gamma$ du $\ell$-groupe cyclique $\Gamma=\Gal(N/K)$ et partons d'une classe invariante $cl^T_S(\A_N) \in \Cl^T_S(N)$, représentée par un $S$-diviseur étranger à $T$. Par hypothèse, le $S$-diviseur $\A_N^{\gamma-1}$ est alors principal, engendré par un élément $T$-infinitésimal $\alpha_N\in\R^T(N)$ défini modulo une $S$-unité $T$-infinitésimale. Sa norme $\varepsilon_K=N_{N/K}(\alpha_N)$ est ainsi une $S$-unité $T$-infinitésimale de $K$ (i.e. un élément de $\E_S^T(K)=\R^T(K)\cap\J^T_S(K)$), définie modulo la norme d'une $S$-unité $T$-infinitésimale de $N$. On obtient par là un morphisme de 
$\Cl^T_S(N)^\Gamma$ vers $\big(\E^T_S(K)\cap N_{N/K}(\R(N))\big)/N_{N/K}(\E_S^T(N))$ dont le noyau est clairement $cl^T_S\big(\D^T_S(N)^\Gamma\big)$ et qui est surjectif en vertu du théorème 90 de Hilbert appliqué aux groupes des $S$-diviseurs étrangers à $T$ (cf. e.g. \cite{J18}, Lem. II.2.32).

\newpage
\section{Module d'Iwasawa $S$-décomposé $T$-ramifié d'un corps surcirculaire}

Montons maintenant la $\Zl$-tour cyclotomique $(K_n)_{n\in\NN}$ au-dessus de $K$ et considérons le corps surcirculaire $K_{\si{\infty}}=\bigcup_{n\in\NN}K_n$.
Désignons par $S$ et $T$ deux ensembles finis disjoints de places finies de $K$; puis, pour toute extension $N$ de $K$ continuons à noter (en l'absence d'ambiguïté) par $S$ et $T$ les ensembles respectifs de places de $N$ au-dessus des précédentes.\par

Faisons choix d'un générateur topologique $\gamma$ du groupe procyclique  $\Gamma=\Gal(K_{\si{\infty}}/K)$ et écrivons $\Lambda=\Zl[[\gamma-1]]$ l'algèbre d'Iwasawa attachée à $\Gamma$.\smallskip

Le groupe de Galois $\,\C^T_S(K_{\si{\infty}})=\Gal(H^T_S(K_{\si{\infty}})/K_{\si{\infty}})$ de la pro-$\ell$-extension abélienne  maximale $H^T_S(K_{\si{\infty}})$ de $K_{\si{\infty}}$ qui est  $S$-décomposée et $T$-ramifiée s'identifie à la limite projective \smallskip

\centerline{$\C^T_S(K_{\si{\infty}})=\varprojlim \,\Cl^T_S(K_n)$}\smallskip

\noindent des pro-$\ell$-groupes de $S$-classes $T$-infinitésimales des corps $K_n$ pour les applications normes.
Il est bien connu que c'est un $\Lambda$-module noethérien, pseudo-isomorphe comme tel à la somme directe\smallskip

\centerline{$\C^T_S(K_{\si{\infty}}) \sim \Lambda^{\rho^{\si{T}}_{\si{S}}} \oplus (\oplus_{i=1}^h \Lambda/P_i\Lambda)$}\smallskip

\noindent d'un $\Lambda$-module libre de dimension finie et d'un nombre fini de quotients monogènes dont les annulateurs respectifs sont engendrés par des polynômes non nuls $P_i$ ordonnés par divisibilité.
Chaque $P_i$ s'écrit comme produit $\ell^{\mu_i}\wP_i$ d'une puissance de $\ell$ et d'un polynôme distingué de degré $\lambda_i$. Le produit $\chi^{\si{T}}_{\si{S}}=\prod P_i$ est le polynôme caractéristique du module d'Iwasawa $\C^T_S(K_{\si{\infty}})$; et les entiers naturels $\rho^{\si{T}}_{\si{S}}$, $\mu^{\si{T}}_{\si{S}}=\sum \mu_i=\nu^\ph_\ell(\chi^{\si{T}}_{\si{S}})$ et $\lambda^{\si{T}}_{\si{S}}=\sum\lambda_i= \deg(\chi^{\si{T}}_{\si{S}})$ sont ses invariants structurels.\smallskip

Le calcul de l'invariant $\rho^{\si{T}}_{\si{S}}$ est effectué dans \cite{JMa} (Th. 9). Il est conjecturé que l'invariant $\mu^{\si{T}}_{\si{S}}$ est toujours nul, ce qui revient à postuler que le polynôme caractéristique $\chi^{\si{T}}_{\si{S}}$ n'est pas divisible par $\ell$. C'est effectivement le cas pour $K$ abélien, lorsque $T$ est vide et que $S$ est l'ensemble $Pl_\ell$ des places au-dessus de $\ell$ (autrement dit lorsque $H^T_S(K_{\si{\infty}})$ est la pro-$\ell$-extension abélienne non ramifiée $\ell$-décomposée maximale de $K_{\si{\infty}}$), en vertu d'un théorème de Ferrero et Washington; et ce résultat vaut encore, sous des conditions moins restrictives sur les ensembles $S$ et $T$ (cf. \cite{JMa}, Th. 13).\smallskip

Enfin, toujours pour $K$ abélien, Greenberg \cite{Grb1} a montré que  $\chi^{\si{\,\emptyset}}_{\si{Pl_\ell}}$ n'est pas non plus divisible par $\omega=\gamma-1$.
Dans ce même contexte abélien, nous allons démontrer la conjecture suivante:

\begin{CoCy}
Soient $\ell$ un nombre premier arbitraire; $K$ un corps de nombres; et $K_{\si{\infty}}=\cup_{n\in\NN}K_n$ sa $\Zl$-extension cyclotomique; soit $Pl_\ell=S \sqcup T$ une partition arbitraire de  l'ensemble $Pl_\ell$ des places au-dessus de $\ell$.\par

Alors le polynôme caractéristique du $\Lambda$-module $\,\C^T_S(K_{\si{\infty}})=\Gal(H^T_S(K_{\si{\infty}})/K_{\si{\infty}})$ attaché à la pro-$\ell$-extension abélienne  maximale $H^T_S(K_{\si{\infty}})$ de $K_{\si{\infty}}$ qui est  $S$-décomposée et $T$-ramifiée n'est pas divisible par $\omega=\gamma-1$.
 En d'autres termes, son sous-module des points fixes  $\,\C^T_S(K_{\si{\infty}})^\Gamma$ est fini:\smallskip
 
 \centerline{$\C^T_S(K_{\si{\infty}})^\Gamma \sim 1$.}
\end{CoCy}

\begin{Prop}
La Conjecture cyclotomique contient celles de Gross-Kuz'min et de Leopoldt.
\end{Prop}

\Preuve Prenons $S=Pl_\ell$ et $T=\emptyset$. Le pro-$\ell$-groupe $\,\C^{\,\emptyset}_{Pl_\ell}(K_{\si{\infty}})$ n'est alors autre que le module de Kuz'min-Tate $\T(K_{\si{\infty}})$ (cf. \cite{J55}, \S2), i.e. le groupe de Galois $ \Gal(K_{\si{\infty}}^{\si{\rm cd}}/K_{\si{\infty}})$ de la pro-$\ell$-extension abélienne maximale de  $K_{\si{\infty}}$ qui est complètement décomposée en toutes les places (car la montée dans la tour cyclotomique $K_{\si{\infty}}/K$ a épuisé toute possibilité d'inertie aux places étrangères à $\ell$). Cela étant, comme $\T(K_{\si{\infty}})$ est un $\Lambda$-module de torsion, la finitude de son sous-groupe invariant $\T(K_{\si{\infty}})^\Gamma$ équivaut à celle de son quotient des genres ${}^\Gamma\T(K_{\si{\infty}})=\T(K_{\si{\infty}})/\T(K_{\si{\infty}})^{\gamma-1}$, lequel est le $\ell$-groupe des classes logarithmiques $\wCl(K)$ dont la conjecture de Gross-Kuz'min pour $K$ postule précisément la finitude (cf. e.g. \cite{J55}, Th. 5). La Conjecture implique donc celle de Gross-Kuz'min.\par

Penons maintenant $S=\emptyset$ et $T=Pl_\ell$. Le pro-$\ell$-groupe $\,\C_\emptyset^{Pl_\ell}(K_{\si{\infty}})$ est alors le  groupe de Galois $ \X(K_{\si{\infty}})=\Gal(K_{\si{\infty}}^{\si{\rm \ell r}}/K_{\si{\infty}})$ de la pro-$\ell$-extension abélienne $\ell$-ramifiée maximale de  $K_{\si{\infty}}$. Et il est bien connu depuis Iwasawa \cite{Iw} que  $\X(K_{\si{\infty}})$ est un $\Lambda$-module qui a pour dimension $\rho^{\si{Pl_\ell}}_{\si{\emptyset}}$ le nombre $c_{\si{K}}$ de places complexes du corps $K$. La finitude du sous-groupe invariant $\X(K_{\si{\infty}})^\Gamma$ s'écrit donc aussi bien ${}^\Gamma\T(K_{\si{\infty}})\simeq \Zl^{c_{\si{K}}}\oplus\T(K)$, pour un certain module fini $\T(K)$. Or, ${}^\Gamma\T(K_{\si{\infty}})$ n'est autre que le groupe de Galois $\Gal(K^{\rm \ell r}/K_{\si{\infty}})$ attaché à la pro-$\ell$-extension abélienne $\ell$-ramifiée maximale de $K$. La finitude de $\T(K)$ exprime donc l'existence d'exactement $(c_{\si{K}}+1)$ $\Zl$-extensions de $K$ linéairement indépendantes, ce qui est précisément la conjecture de Leopoldt (cf. e.g. \cite{J55}, Th. 12). 

\newpage
\section{Équivalence avec les conjectures de Leopoldt et de Gross-Kuz'min}

Regardons d'abord le cas totalement réel (déjà étudié dans \cite{J43} \S3 via la théorie des genres):

\begin{Th}
Pour $K$ totalement réel, la Conjecture cyclotomique résulte de celle de Leopoldt.
\end{Th}

\Preuve La conjecture de Leopoldt affirme qu'un corps totalement réel $K$ possède une unique $\Zl$-extension, autrement dit que sa pro-$\ell$-extension abélienne $\ell$-ramifiée maximale $K^{\rm \ell r}$ est de degré fini sur la $\Zl$-extension cyclotomique $K^{\rm c}=K_{\si{\infty}}$.
Considérons alors le quotient des genres:\smallskip

\centerline{${}^\Gamma\C^T_S(K_{\si{\infty}})=\C^T_S(K_{\si{\infty}}) / \big(\C^T_S(K_{\si{\infty}})\big)^{\gamma-1}$}\smallskip

\noindent du groupe $\,\C^T_S(K_{\si{\infty}})$. Par construction, c'est le groupe de Galois $\Gal(H^T_S(K_{\si{\infty}}/K)/K_{\si{\infty}})$ de la plus grande sous-extension $H^T_S(K_{\si{\infty}}/K)$ de $H^T_S(K_{\si{\infty}})$ qui est abélienne sur $K$. Or, celle-ci est $T$-ramifiée sur $K_{\si{\infty}}$, donc $\ell$-ramifiée sur $K$, i.e. contenue dans $K^{\rm \ell r}$. Ainsi ${}^\Gamma\C^T_S(K_{\si{\infty}})$ est fini; de même
$\,\C^T_S(K_{\si{\infty}})^\Gamma$.\medskip

Supposons maintenant que $K$ soit une extension quadratique totalement imaginaire d'un sous-corps totalement réel $K^+$. Notons $\tau$ la conjugaison complexe et $\Delta=\Gal(K/K^+)=\{1,\tau\}$.
\begin{itemize}
\item Si $\ell$ est impair, les idempotents orthogonaux $e_\pm= \frac{1}{2}(1 \pm \tau)$ permettent de décomposer canoniquement chaque $\Zl[\Delta]$-module $M$ comme somme directe de ses composantes réelle $M^+=M^{e_{\si{+}}}$ et imaginaire $M^-=M^{e_{\si{-}}}$. On a ainsi: $M=M^+\oplus M^-$.
\item Si $\ell$ vaut $2$ et que $M$ est $\Z2$-noethérien, confondre le noyau de $(1\pm\tau)$ avec l'image de $(1\mp\tau)$ donne lieu à une erreur finie.
On peut donc continuer à définir composantes réelle et imaginaire de $M$ comme image et noyau respectifs de $(1+\tau)$ et écrire à un fini près:
\centerline{$M \sim M^+\oplus M^-$.}
\end{itemize}

Avec ces conventions,  le {\em Théorème principal} de cette note s'énonce comme suit:

\begin{Th}\label{ThP} 
Soient $K$ une extension quadratique totalement imaginaire d'un corps totalement réel $K^+$ et $S \sqcup T$ une partition de  l'ensemble $Pl_\ell$ des places de $K^+$ au-dessus de $\ell$.
Décomposons  $\,\C^T_S(K_{\si{\infty}})^\Gamma$ (à un fini près pour $\ell=2)$ en ses composantes réelle et imaginaire. Alors:
\begin{itemize}
\item[(i)] La finitude de $(\,\C^T_S(K_{\si{\infty}})^\Gamma)^+$ résulte de la conjecture de Leopoldt pour $K$.
\item[(ii)] La finitude de $(\,\C^T_S(K_{\si{\infty}})^\Gamma)^-$ résulte de la conjecture de Gross-Kuz'min pour $K$.
\end{itemize}
En d'autres termes, la Conjecture cyclotomique pour le corps $K$, i.e. la finitude de $\,\C^T_S(K_{\si{\infty}})^\Gamma$, est équivalente à la réunion des conjectures de Leopoldt et de Gross-Kuz'min pour ce même corps.
\end{Th}

\Preuve Examinons successivement les deux assertions:

$(i)$ provient directement du cas totalement réel traité ci-dessus, la composante réelle du groupe  $\,\C^T_S(K_{\si{\infty}})^\Gamma$ correspondant (à un fini près pour $\ell=2$) au même groupe  $\,\C^T_S(K^+_{\si{\infty}})^\Gamma$ pour le corps $K^+$.

$(ii)$ résulte de la généralisation suivante d'un théorème de Kuz'min (cf. \cite{Kuz}, Prop. 7.5 et \cite{J55}, Th. 17) et du fait que, sous la conjecture de Gross-Kuz'min dans $K$, la composante imaginaire du $\ell$-groupe des unités logarithmiques $\,\wE(K)$ se réduit au $\ell$-groupe des racines de l'unité (cf. \cite{J28}, \S 3):

\begin{Th}
Étant donnés deux ensembles finis disjoints $S$ et $T$ de places d'un corps de nombres $K$, si $S \cup T$ contient l'ensemble $Pl_\ell$ des places au-dessus de $\ell$, on a un isomorphisme naturel:\smallskip

\centerline{$\,\C^T_S(K_{\si{\infty}})^\Gamma \simeq \big(\E^T_S(K) \cap \,\wE(K)\big) / \E^T_S(K)^\nu$,}\smallskip

\noindent où $\,\E^T_S(K) \cap \,\wE(K)$ est l'intersection du groupe  $\,\E^T_S(K) $ des $S$-unités $T$-infinitésimales de $K$ avec le groupe $\wE(K)$ des unités logarithmiques; et  $\,\E^T_S(K)^\nu=\bigcap_{n\in\NN}N_{K_n/K}(\E^T_S(K_n))$ est le sous-groupe (dit universel) formé des éléments qui sont normes de $S$-unités $T$-infinitésimales à tous les étages.
\end{Th}

\Preuve Partons de l'isomorphisme donné par la Proposition \ref{Prop} et passons à la limite projective pour les applications normes.
Côté diviseurs à gauche, la réunion $S\cup T$ contenant par hypothèse l'ensemble des places ramifiées dans la tour $K_{\si{\infty}}/K$, les $S$-diviseurs ambiges et étrangers à $T$ de chacun des corps $K_n$ se réduisent aux seul diviseurs provenant de $K$. On a donc tout simplement:\smallskip

\centerline{$\varprojlim \big(\Cl^T_S(K_n)^\Gamma/cl^T_S\big(\D^T_S(K_n)^\Gamma\big)\big)   = \varprojlim \,\Cl^T_S(K_n)^\Gamma  = \, \C^T_S(K_{\si{\infty}})^\Gamma$.}\smallskip

\noindent Côté unités à droite, l'intersection $\bigcap_{n\in\NN}N_{K_n/K}(\R(K_n))$ caractérise précisément le sous-groupe $\,\wE(K)$ des unités logarithmiques dans $\R(K)$ (cf. \cite{J28}, Prop. 3.2 ou \cite{J55}, \S1). D'où le résultat.

\newpage
\section{Conséquences du Théorème principal}

Il résulte du Théorème \ref{ThP} que la Conjecture cyclotomique vaut pour tous les corps $K$ réels ou à conjugaison complexe (avec, dans ce cas, $S$ et $T$ stables par la conjugaison) et les premiers $\ell$ pour lesquels sont simultanément satisfaites les conjectures de Leopoldt et de Gross-Kuz'min. Ainsi:

\begin{Cor}
La Conjecture cyclotomique est satisfaite en particulier par:
\begin{itemize}
\item les corps de nombres abéliens $K$, pour n'importe quel premier $\ell$;
\item les extensions quadratiques totalement imaginaires $K$ d'un corps totalement réel qui contien\-nent les racines $2\ell$-ièmes de l'unité et sont $\ell$-logarithmiquement principales (i.e. dont le $\ell$-groupe des classes logarithmiques $\,\wCl(K)$ est trivial), pour ce même premier $\ell$.
\end{itemize}

Dans chacun de ces deux cas, la Conjecture est alors vérifiée par chacun des étages finis $K_n$ de la $\Zl$-tour cyclotomique $K_{\si{\infty}}$; et le polynôme caractéristique $\chi_{\si{S}}^{\si{T}}$ du $\Lambda$-module d'Iwasawa $\,\C^T_S(K_{\si{\infty}})$ n'est donc divisible par aucun des polynômes cyclotomiques $\omega_n=(\gamma^{\ell^{\si{n}}}-1)/(\gamma^{\ell^{\si{n-1}}}-1)$.
\end{Cor}

\Preuve Le premier cas résulte du fait que les corps abéliens satisfont les conjectures de Leopoldt et de Gross-Kuz'min pour tous les nombres premiers $\ell$ en vertu du théorème d'indépendance de logarithmes de nombres algébriques de Baker-Brumer. Le second cas provient du fait que les corps $\ell$-logarithmiquement principaux vérifient banalement la conjecture de Gross-Kuz'min (qui postule la finitude du groupe $\,\wCl(K)$ pour le nombre premier $\ell$); et qu'ils vérifient en outre celle de Leopoldt s'ils contiennent les racines $2\ell$-ièmes de l'unité (cf. e.g. \cite{J17} ou \cite{J31}).\par
Enfin, si $K$ est abélien, tous les $K_n$ le sont. Et s'il est à conjugaison complexe et $\ell$-logarithmique\-ment principal, tous les $K_n$ le sont aussi; d'où le résultat annoncé.
\medskip

Soient maintenant $S$ et $T$ deux ensembles finis disjoints arbitraires de places du corps  $K$. Notons $S_\ell=S \cap Pl_\ell$ la partie {\em sauvage} de $S$ et $S_\circ=S\setminus S_\ell$ sa partie {\em modérée}. Écrivons de même $T=T_\ell \cup T_\circ$.
Notons toujours $\,\C^T_S(K_{\si{\infty}})=\Gal(H^T_S(K_{\si{\infty}})/K_{\si{\infty}})$ le groupe de Galois de la pro-$\ell$-extension abélienne  maximale $H^T_S(K_{\si{\infty}})$ de $K_{\si{\infty}}$ qui est  $S$-décomposée et $T$-ramifiée. 

\begin{Sco}\label{Sco} 
Avec ces conventions, il vient:
\begin{itemize}
\item[(i)] Le groupe $\,\C^T_S(K_{\si{\infty}})$ est indépendant de $S_\circ$; autrement dit, on a: $\,\C^T_S(K_{\si{\infty}})=\,\C^T_{S_{\si{\ell}}}(K_{\si{\infty}})$.
\item[(ii)] Les groupes $\,\C^T_S(K_{\si{\infty}})$ et $\,\C^{T_{\si{\ell}}}_S(K_{\si{\infty}})$ ont même dimension comme $\Lambda$-modules et leurs quotients des genres sont pseudo-isomorphes: ${}^\Gamma\C^T_S(K_{\si{\infty}}) \sim {}^\Gamma\C^{T_{\si{\ell}}}_S(K_{\si{\infty}})$.
\end{itemize}
Il en résulte un pseudo-isomorphisme entre sous-groupes invariants: $\,\C^T_S(K_{\si{\infty}})^\Gamma \sim \,\C^{T_{\si{\ell}}}_{S_{\si{\ell}}}(K_{\si{\infty}})^\Gamma$.
Et la Conjecture cyclotomique revient donc à postuler la finitude de $\,\C^{T}_{S}(K_{\si{\infty}})^\Gamma$ pour $Pl_\ell \subset S_\ell \cup T_\ell$.
\end{Sco}

\Preuve Les places modérées sont presque totalement inertes dans la $\Zl$-tour cyclotomique $K_{\si{\infty}}/K$.\par

$(i)$ De ce fait, la montée dans la tour ayant épuisé toute possibilité d'inertie, une place modérée non-ramifiée  au-dessus de $K_{\si{\infty}}$ est donc complètement décomposée. Il suit de là que $H_S^T(K_{\si{\infty}})$ coïncide avec $H_{S_{\si{\ell}}}^T(K_{\si{\infty}})$ et $\,\C^T_S(K_{\si{\infty}})=\Gal(H_{S_{\si{\ell}}}^T(K_{\si{\infty}})/K_{\si{\infty}})$ avec $\,\C^T_{S_{\si{\ell}}}(K_{\si{\infty}})=\Gal(H_{S_{\si{\ell}}}^T(K_{\si{\infty}})/K_{\si{\infty}})$.\par

$(ii)$ Le cas de $T$ est plus compliqué: 

-- D'un côté, puisque les places de $T_\circ$ sont finiment décomposées dans la tour cyclotomique, les groupes finis d'unités semi-locales $\,\U_{T_\circ\!}(K_n)=\prod_{\p_n\in T_\circ}\,\U_{K_{\p_{\si{n}}}}$ sont de rang borné. Pour toute pro-$\ell$-extension abélienne donnée $H_{\si{\infty}}$ de $K_{\si{\infty}}$ le sous-groupe de $\Gal(H_{\si{\infty}}/K_{\si{\infty}})$ engendré par les sous-groupes d'inertie attachés aux places de $T_\circ$ est ainsi un $\Zl$-module de type fini et donc un $\Lambda$-module de torsion. En particulier  $\,\C^T_S(K_{\si{\infty}})$ et $\,\C^{T_{\si{\ell}}}_S(K_{\si{\infty}})$ ont même dimension comme $\Lambda$-modules.

-- D'un autre côté, le groupe semi-local $\,\U_{T_\circ\!}(K)=\prod_{\p\in T_\circ}\,\U_{K_\p}$ étant fini, son  image dans le groupe de Galois ${}^\Gamma\C^T_S(K_{\si{\infty}})=\Gal(H^T_S(K_{\si{\infty}}/K)/K_{\si{\infty}})$ attaché à la plus grande sous-extension $H^T_S(K_{\si{\infty}}/K)$ de $H^T_S(K_{\si{\infty}})$ qui est abélienne sur $K$ l'est aussi; et les quotients  ${}^\Gamma\C^T_S(K_{\si{\infty}})$ et $ {}^\Gamma\C^{T_{\si{\ell}}}_S(K_{\si{\infty}})$ sont donc pseudo-isomorphes.

-- Réunissant ces résultats, on conclut que les sous-groupes invariants ${}^\Gamma\C^T_S(K_{\si{\infty}})$ et ${}^\Gamma\C^{T_{\si{\ell}}}_S(K_{\si{\infty}})$ sont eux-mêmes pseudo-isomorphes.\medskip

\Remarque Sauf à compliquer la Conjecture (à l'instar de \cite{LS}), la condition $Pl_\ell \subset S_\ell \cup T_\ell$ reste, en revanche, incontournable: cf. e.g. \cite{J18}, Th. IV.2.9, pour un contre-exemple.

\newpage

\noindent {\em Commentaires bibliographiques}\smallskip

 Les $\ell$-groupes de $S$-classes $T$-infinitésimales sont étudiés dans \cite{J18} (Ch. II, \S 2). On y trouve notamment la suite exacte des classes ambiges évoquée dans la section 1.\par

Les $\ell$-groupes de classes logarithmiques ont été introduits dans \cite{J28}. Leur calcul est maintenant implanté dans {\sc pari} (cf. \cite{BJ}).\par
 
Les principaux résultats de la Théorie $\ell$-adique du corps de classes introduite dans \cite{J18} sont présentés dans \cite{J31}. On peut aussi se reporter au livre de Gras \cite{Gra2}.\par

Enfin, le calcul des invariants d'Iwasawa attachés aux $\ell$-groupes de $S$-classes $T$-infinitésimales est développé dans \cite{J43,JMa,JMP} en liaison avec les identités de dualité de Gras.

\bigskip
\def\refname{\normalsize{\sc  Références}}

\medskip\noindent
{\small
\begin{tabular}{l}
Institut de Mathématiques de Bordeaux \\
Université de {\sc Bordeaux} \& CNRS \\
351 cours de la libération\\
F-33405 {\sc Talence} Cedex\\
courriel : Jean-Francois.Jaulent@math.u-bordeaux.fr\\
\url{https://www.math.u-bordeaux.fr/~jjaulent/}
\end{tabular}
}

 \end{document}